\documentstyle[12pt,amscd,leqno]{amsart}
\newtheorem{Theorem}{Theorem}[section]
\newtheorem{Lemma}[Theorem]{Lemma}

\newtheorem{Proposition}[Theorem]{Proposition}
\newtheorem{Example}[Theorem]{Example}
\newtheorem{Remark}[Theorem]{Remark}

\def\Ass{\operatorname{Ass}}
\def\Det{\operatorname{Det}}
\def\deg{\operatorname{deg}}

\def\Nset{{\mathbb N}}
\def\Zset{{\mathbb Z}}
\def\Rset{{\mathbb R}}
\def\mm{{\mathfrak m}}
\def\pfrak{{\mathfrak p}}
\def\Ccal{{\mathcal C}}
\def\Ical{{\mathcal I}}
\def\Rcal{{\mathcal R}}
\def\Acal{{\mathcal A}}
\def\ada{{\mathbf a}}
\def\alphda{\underline{\alpha}}
\def\betda{\underline{\beta}}
\def\bda{{\mathbf b}}
\def\cda{{\mathbf c}}

\def\xda{{\mathbf x}}
\def\uda{{\mathbf u}}
\def\vda{{\mathbf v}}
\def\wda{{\mathbf w}}
\def\tda{{\mathbf t}}

\title{ Stability of associated primes of monomial ideals}
\author{  L\^e Tu\^an Hoa }
\address{ Institute of Mathematics Hanoi, 18 Hoang Quoc Viet Road, 10307 Hanoi, Vietnam  }
\email{lthoa@@math.ac.vn}
\begin{document}
\thanks{The  author was  supported  in part by the National Basic Research
Program (Vietnam). } 
\begin{abstract}  Let $I$ be a monomial ideal of a polynomial ring $R$. In
this paper we determine a number $B$ such that $\Ass (I^n/I^{n+1}) = \Ass (I^{B}/I^{B+1})$
for all $n\geq B$.
\end{abstract}

 \maketitle

\rightline{\it Dedicated to Professor Do Long Van}

\rightline{\it on the occasion of his 65-th birthday.}

\section*{Introduction}\smallskip

Let $I$ be an ideal of a Noetherian ring $R$.  It is a well-known result of M.
Brodmann \cite{B} that the sequences $\{ \Ass (R/I^n)\}_{n\geq 1}$ and $\{
\Ass (I^n/I^{n+1})\}_{n\geq 1}$ stabilize for large $n$. That is, there are
positive numbers $A$ and $B$ such that $\Ass (R/I^n) = \Ass (R/I^A)$ for all
$n\geq A$ and $ \Ass (I^n/I^{n+1}) = \Ass (I^B/I^{B+1})$ for all $n\geq B$.
Very little is known about the numbers $A$ and $B$. One of the difficulties in
estimating these numbers is that neither of the above sequences is monotonic; see
\cite{ME} and also \cite{HH} for monomial ideals. In an earlier
paper of S. McAdam and P. Eakin \cite{ME} and a recent paper of R. Sharp
\cite{Sh} there are some information about the behavior of these sequences.
Moreover, for  specific prime ideals $\pfrak$ one can  decide in
terms of the Castelnuovo-Mumford regularity of the associated graded
ring of $I$ when $\pfrak$ belongs to $\Ass (R/I^n)$ (see \cite{Sh}, Theorem
2.10). For a very restricted class of ideals the numbers $A$ and $B$ can be
 rather small (see \cite{M}).

The aim of this paper is to find an explicit value for $A$ and $B$ for a
monomial ideal $I$ in a polynomial ring $R =K[t_1,...,t_r]$  over a field.  A special case was studied
in \cite {CMS}, when $I$ is generated by products of
two different variables. Such an ideal is associated to a graph. The result looks
nice: the number $A$ can be taken as the number of variables (see \cite{CMS},
Proposition 4.2, Lemma 3.1 and Corollary 2.2). However the  approach  of \cite{CMS} cannot be applied for arbitrary monomial ideals.

 It is interesting to note that in our situation we always have $A=B$, since $\Ass (R/I^n) = \Ass
(I^{n-1}/I^n) $  (see \cite{T}, Proposition 5).  In this paper, it is more convenient for us to work with
$\Ass (I^n/I^{n+1})$ (and hence with number $B$). Let $\mm = (t_1,...,t_r)$.  Then one can reduce
the problem of finding $B$ to finding a number $B'$ such that $\mm \in \Ass (I^n/I^{n+1}) $ for all
$n\geq B'$ or  $\mm \not\in \Ass (I^n/I^{n+1} ) $ for all $n\geq B'$ (see  Lemma \ref{B1}).  From
this
observation we have to study the vanishing (or non-vanishing) of the local
cohomology module $H^0_\mm(I^n/I^{n+1} )$. The main technique to
do that is to describe these sets as graded components of certain modules over
toric rings raised from systems of linear constraints. Then we have to bound the
degrees of generators of these modules, and also to bound certain invariants related to
the Catelnuovo-Mumford regularity. The numbers $B$ found in Theorem
\ref{D10} depends on the number of
variables $r$, the number of generators $s$ and the maximal degree $d$ of
generators of $I$. This number is very big. However there are examples
 showing  that such a number $B$ should also involve $d$ and $r$ (see Examples
\ref{E11} and \ref{E12}).

The paper is divided into two sections. The first one is of preparatory
character. There we will give a bound for the degrees of generators of a
module raised from integer solutions of a system of linear constraints.
Section 2 is  devoted to determining the number $B$. First we will find a
number from which the sequence $\{ \Ass (I^n/I^{n+1})\}_{n\geq 1}$ is
decreasing (see Proposition \ref{D5}). Then we will have to bound a number
related to the Castelnuovo-Mumford regularity of the associated graded ring of
$I$ (Proposition \ref{D9}) in order to use a result of \cite{ME} on the
increasing property of this sequence. The main result of the paper  is Theorem \ref{D10}.
This   section will be  ended with two examples which  show how big
$B$ should be.

I would like to end this introduction with the remark that by a different
method, T. N. Trung \cite{T} is able to solve similar problems for the
integral closures of powers of a monomial ideal.

\section{Integer solutions of linear constraints}\smallskip

Let $S$ be  the set of integer solutions of the following
system of linear constraints
\begin{eqnarray}\label{EA1}
\begin{cases}
a_{i1} x_1 + \cdots + a_{ie}x_e \geq 0, \ (i=1,...,s),\\
x_1\geq 0,..., x_e \geq 0,
\end{cases}
\end{eqnarray}
where $a_{ij}\in \Zset$. It is a fundamental fact in integer programming that
the semigroup ring $K[S]$   is a finitely generated subring of
$K[x_1,...,x_e]$. An algebraic proof can be found in \cite{S}, Section 1.3.
What we need is an ``effective" version of this result. To this end we will
consider an element of $S$ as a point in the space $\Rset^e$. For a vector
$\vda = (v_1,...,v_e) \in \Rset^e$, put
$$ \| \vda \| = \sqrt{v_1^2 + \cdots + v_e^2} \ \ \text{and} \ \ \| \vda \|_* = \max
\{ |v_1|,...,|v_e|\}.$$

The proof of the following lemma and Lemma \ref{A2} is similar to that of
\cite{Sch}, Theorem 17.1. For  convenience of the readers we give here the
detail.

\begin{Lemma}\label{A1} Let $\ada_j = (a_{1j},...,a_{sj})^T \in \Zset^s$
denote the coefficient column of $x_i$ in (\ref{EA1}). Assume that $\|\ada_1
\| \geq \cdots \geq \| \ada_e \| >0$. Then $K[S]$ is generated by monomials
$\xda^\vda := x_1^{v_1}\cdots x_e^{v_e}$ such that
$$\|\vda \|_* < e\|\ada_1\| \cdots \|\ada_{e-1}\| \leq e\|\ada_1\| \cdots \|\ada_e \|.$$
\end{Lemma}

 \begin{pf}
 Let $\Ccal$ be the set of all real solutions of (\ref{EA1}). It is a polyhedral convex
 set in $\Rset^e$. By Minkovski's Theorem (see \cite{Sch}, Corollary 7.1a), one
 can write
 $$ \Ccal = \Rset_+ \uda_1 + \cdots + \Rset_+ \uda_k,$$
 where $\uda_1,..., \uda_k \in \Ccal$ and $\Rset_+$ is the set of nonnegative numbers.
 Here we choose $k$ the smallest possible. Then $\Rset_+ \uda_1 , ..., \Rset_+ \uda_k$
 are extreme rays.
 Each  extreme ray is an intersection of $e-1$ independent hyperplanes
 appeared in (\ref{EA1}). Hence, we may without loss of generality assume that
  $\uda_p,\ 1\leq p \leq k$, is a nonzero solution of a linear subsystem of
  the type
\begin{eqnarray}\label{EA2}
\begin{cases}
b_{i1}x_1+ \cdots + b_{iq}x_q = - b_{i, q+1}x_{q+1},  \ (i=1,...,q),\\
x_{q+2} = \cdots = x_e =   0,
\end{cases}
\end{eqnarray}
where $q \leq \min\{ e-1, s\}$, the matrix on the left hand is invertible, and
each column vector $\bda_j$ is a subvector of $\ada_j$. By Cramer's rule we
may choose $\uda_p$ the integer solution:
$$\uda_p = (D_1,...,D_q,D_{q+1},0,...,0),$$
where $D_1,...,D_{q+1}$ are determinants of the linear system consisting of
the first $q$ equations of (\ref{EA2}). Note that if $\cda_1,...,\cda_q \in
\Rset^q$ are column vectors, then
\begin{eqnarray}\label{EA3b}
\Det (\cda_1,...,\cda_q) \leq \|\cda_1\|\cdots \|\cda_q\|.
\end{eqnarray}
Hence
$$\|\uda_p\|_* \leq \max_i \|\bda_1\|\cdots \|\bda_{i-1}\|\|\bda_{i+1}\|
\cdots \|\bda_{q+1}\|\leq \|\ada_1\|\cdots \|\ada_q\| \leq \|\ada_1\|\cdots
\|\ada_{e-1}\|.$$ From now on we assume that all elements $\uda_1,...,\uda_k$
are integer points chosen in the above way. In particular they belong to $S$.
Let $\vda\in S$ be an arbitrary element. Since $\vda \in \Ccal$, by
Caratheodory's Theorem (see \cite{Sch}, Corollary 7.1i), one can find $\{
i_1,...,i_q \} \subseteq \{ 1,...,k\},\ q \leq e$ and numbers $ \alpha_{i_1} ,..., \alpha_{i_q}  \geq 0$, such that
$$\vda = \alpha_{i_1}\uda_{i_1} + \cdots + \alpha_{i_q}\uda_{i_q}.$$
 For a real number $\alpha$, let $[\alpha]$ denote
the largest integer not exceeding $\alpha$. Let
$$\uda =
[\alpha_{i_1}]\uda_{i_1} + \cdots + [\alpha_{i_q}]\uda_{i_q},$$ and
$$\wda =
(\alpha_{i_1} - [\alpha_{i_1}])\uda_{i_1} + \cdots + (\alpha_{i_q} -
[\alpha_{i_q}])\uda_{i_q}.$$ We have $\wda = \vda - \uda \in \Nset^e$. However
$\wda \in \Ccal$. Hence $\wda \in \Ccal \cap \Nset^e = S.$ Since $\vda = \uda
+ \wda$, this means that the following set generates $S$:
$$\{ \uda_1,...,\uda_k\} \cup  \{\alpha_{i_1}\uda_{i_1} + \cdots +
\alpha_{i_q}\uda_{i_q} \in \Nset^e | \ q\leq e,\ 0 \leq \alpha_{i_j} <1,\
1\leq i_1,...,i_q\leq k \}.$$ For each vector $\vda = \alpha_{i_1}\uda_{i_1} +
\cdots + \alpha_{i_q}\uda_{i_q}$ in the second subset of the above union we
have
$$\| \vda\|_* < q \max_{j=1,...,k} \|\uda_j\|_* \leq e \max_{j=1,...,k}
\|\uda_j\|_* \leq e \|\ada_1\|\cdots \|\ada_{e-1}\|.$$
Hence the assertion
holds true.
\end{pf}

The following simple example shows that the above result is essentially
optimal.

\begin{Example} {\rm
Consider the system of constraints
$$\left\{ \begin{array}{rl}
dx_1 - x_2 & \geq 0, \\
\cdots & \\
dx_{e-1} - x_e & \geq 0, \\
x_1\geq 0,...,x_e & \geq 0.
\end{array}\right.$$
The corresponding polyhedral convex set has an extreme ray $\Rset_+\uda$,
where $\uda = (1,d,...,d^{e-1})$. Clearly, $\uda$ is a minimal generator of
$S$.}
\end{Example}

We now consider the set $E$ of integer solutions of the
following system of linear constraints:
\begin{eqnarray}\label{EA4}
\begin{cases}
a_{i1} x_1 + \cdots + a_{ie}x_e \geq b_i, \ (i=1,...,s),\\
x_1\geq 0,..., x_e \geq 0,
\end{cases}
\end{eqnarray}
where $a_{ij}, b_i \in \Zset$. Since $S + E \subseteq E$, $K[E]$ is a module
over $K[S]$. For simplicity, sometimes we also say that $E$ is a $S$-module.

\begin{Lemma}\label{A2} Keep the notation of Lemma \ref{A1}. Let $\bda = (b_1,...,b_s)^T \in
\Zset^s$. Then the module $K[E]$ is generated over $K[S]$ by monomials
$\xda^\vda $ such that
$$\|\vda \|_* < (e + \|\bda \|)\|\ada_1\| \cdots \|\ada_e\|.$$
\end{Lemma}

\begin{pf} Let $\Ccal'$ be the set of all real solutions of (\ref{EA4}). Then
$\Ccal'$ is also a polyhedral convex set. By Minkovski's Theorem one can write
$$\Ccal' = \{ \lambda_1\uda_1+\cdots +\lambda_k\uda_k + \mu_1 \vda_1+ \cdots +
\mu_l\vda_l | \ \lambda_i,\mu_j \geq 0,\ \sum \mu_j = 1\},$$ where
$\uda_1,...,\uda_k$ are defined in the proof of the previous lemma, and
$\vda_1,...,\vda_l$ are extreme points. These extreme points are solutions of
$e$ independent affine hyperplanes appeared in (\ref{EA4}). By a similar
argument to the proof of Lemma \ref{A1} we get that
$$\| \vda_j \|_* \leq \|\bda\|\cdot \|\ada_1\| \cdots \|\ada_{e-1}\|,$$
and that the set
$$\begin{array}{ll}
\{ \lambda_{i_1}\uda_{i_1}+\cdots +\lambda_{i_q}\uda_{i_q} + \mu_1 \vda_1+
\cdots + \mu_l\vda_l \in \Nset^e | & q\leq e,\ 1\leq i_1< \cdots < i_q \leq
k,\\ & 0\leq \lambda_{i_j} < 1,\ \mu_j\geq 0,\ \sum \mu_j = 1\}
\end{array}$$
generates the module $E$ over $S$. All these elements have the $*$-norms less
than
$$e \max_i \|\uda_i\|_* + \max_j \|\vda_j\|_* \leq (e + \|\bda \|)\|\ada_1\|
 \cdots \|\ada_e\|,$$
which proves the assertion. \end{pf}

\noindent {\bf Remark}. In the sequel, by abuse of terminology, if
$$\varphi (\xda) = a_1x_1+ \cdots + a_ex_e,$$
is a linear functional, then we say that $\varphi(\xda) \geq 0$ is a
homogeneous linear constraint, while $\varphi(\xda) \geq b$ is a linear
constraint.

\section{Stability of  $ \Ass (I^n/I^{n+1})$ }\smallskip

We always assume that  $I$ is a non-zero  monomial ideal of a polynomial ring $R=
K[t_1,...,t_r]$. If $r\geq 2$, then for a positive integer $j\leq r$ and $\ada=
(a_1,...,a_r) \in \Rset^r$ we set
$$ \ada[j] = (a_1,...,a_{j-1}, a_{j+1},...,a_r).$$
Thus the monomial $\tda^{\ada[j]}$ is obtained from $\tda^{\ada}$ by setting
$t_j=1$. Let $I[j]$ be the ideal generated by all monomials $\tda^{\ada[j]}$
such that $\tda^{\ada}\in I$. Note that
$\tda^{\ada_1[j]},...,\tda^{\ada_s[j]}$ generate $I[j]$ provided $\{
\tda^{\ada_1},...,\tda^{\ada_s} \}$ is a generating system of $I$. Hence for
all $n$ we have
$$I^n[j] = I[j]^n.$$
The following observation is simple but useful. It comes from the fact that
any associated prime of a monomial ideal is generated by a subset of
variables.

\begin{Lemma}\label{B1} Let $\mm = (t_1,...,t_r)$ and $r\geq 2$. Then for all
$n\geq 1$ we have
$$ \Ass(I^n/I^{n+1}) \setminus \{\mm \} = \cup_{i=1}^r \Ass(I[i]^n/I[i]^{n+1}).$$
\end{Lemma}

 \begin{pf} It  immediately follows from \cite{T}, Lemma 11 and Proposition 4.  Another way is to  
modify the proof of Lemma 11 in  \cite{T}.
\end{pf}

Using this lemma, by the  induction on the number of variables, it is clear that in order to study the
stability of $\Ass(I^n/I^{n+1})$  we have to find a number $n_0$ such that $
\mm \in \Ass(I^n/I^{n+1})$ for all $n\geq n_0$, or vice-verse, $ \mm \not\in
\Ass(I^n/I^{n+1})$ for all $n\geq n_0$. Note that $ \mm \in \Ass(I^n/I^{n+1})$
if and only if the local cohomology module $H^0_{\mm}(I^n/I^{n+1}) \neq 0$.
Let
$$G = \oplus_{n\geq 0}I^n/I^{n+1}$$
denote the associated graded ring of $I$. Then $H^0_{\mm G}(G)$  is a graded
$G$-module. Moreover, as a submodule of $G$,  it is a finitely generated
module. We have

\begin{Lemma}\label{D1} For $r\geq 2$,   $H^0_{\mm G}(G)_{n-1} \cong 
H^0_{\mm}(I^{n-1}/I^n) \cong
 \frac{I^{n-1} \cap I[1]^n \cap \cdots \cap
I[r]^n}{I^n}$.
\end{Lemma}

\begin{pf} The first isomorphism  is well-known (see, e.g., \cite{CHT}, Lemma 2.1 for a proof),
 while the second one follows from the fact
$$I^n : (x_1,...,x_r)^\infty = \cap_{i=1}^r (I^n: x_i^\infty) =
\cap_{i=1}^rI[i]^n.$$ Here we denote $I : J^\infty = \cup_{m=1}^\infty I: J^m$.
\end{pf}

The first isomorphism of the above lemma allows us to study
$H^0_{\mm}(I^n/I^{n+1}),\ n\geq 0$, in the total. Our preliminary task is to
bound the degree of generators of the module $H^0_{\mm G}(G)$. Let
$$J = I[1]^n \cap \cdots\cap I[r]^n.$$
 We will try to associate the set of monomials in $J\cap I^{n-1}$ to the set
 of integer solutions of a system of linear constrains, so that we can use the
 results of Section 1. Our technique is based on the following remarks
 which will be used several times. Note that this technique was used in
 Section 7 of \cite{F}.

\begin{Remark}\label{D2b}{\rm \begin{itemize} \item[(i)] An intersection of monomial
 ideals and a quotient of two monomial ideals are again monomial ideals.

\item[(ii)] A monomial ideal is entirely defined by the set of its monomials.
If $I_1 \subset I_2$ are monomial ideals, then the number of monomials in
$I_2\setminus I_1$ is equal to the dimension of the $K$-vector space
$I_2/I_1$.

\item[(iii)] Assume that the  monomials $\tda^{\ada_1},...,
\tda^{\ada_s}$ generate the  ideal $I$. Then a monomial $\tda^{\bda}\in I^n$ if
and only if there are nonnegative integers $\alpha_1,...,\alpha_{s-1}$,  such
that $n\geq \alpha_1+\cdots + \alpha_{s-1}$ and $\tda^{\bda}$ is divisible by
$$(\tda^{\ada_1})^{\alpha_1}\cdots (\tda^{\ada_{s-1}})^{\alpha_{s-1}}
(\tda^{\ada_s})^{n- \alpha_1-\cdots - \alpha_{s-1}}.$$
This is equivalent to
$$b_j \geq a_{1j}\alpha_1 + \cdots + a_{(s-1)j}\alpha_{s-1} + a_{sj}(n- \alpha_1-\cdots -
\alpha_{s-1}),$$ for all $j=1,...,r$, where $\ada_i = (a_{i1},...,a_{ir})$.
\end{itemize}}
\end{Remark}

From now on assume that $I$ is minimally generated by the monomials
$\tda^{\ada_1},..., \tda^{\ada_s}$. Note that if $I$ is generated by powers of
variables, i. e.   $I = (t_{i_1}^{a_1},..., t_{i_p}^{a_p})$, then
$$\Ass(I^n/I^{n+1})  = \{ (t_{i_1},...,t_{i_p}) \}$$
for all $n>0$. Therefore, in the whole paper we may assume that
\vskip0.3cm

($\sharp$) $\ada_s$ contains at least two non-zero components.
\vskip0.3cm

\noindent This will simplify our calculation.

Consider the following system of linear constraints
\begin{eqnarray}\label{ED1}
\begin{cases}
y_j \geq a_{1j}x_1 + \cdots + a_{(s-1)j}x_{s-1} + a_{sj}(z-x_1-\cdots -
x_{s-1}-1), \ (j=1,...,r),\\
z \geq x_1+ \cdots + x_{s-1}+ 1, \\
y_j \geq a_{1j}x_{i1} + \cdots + a_{(s-1)j}x_{i(s-1)} + a_{sj}(z-x_{i1}-\cdots
- x_{i(s-1)}), \\ (i,j=1,...,r;\ j\neq i),\\
z \geq x_{i1}+ \cdots + x_{i(s-1)}, \ (i=1,...,r), \\
z \geq 0;\ y_1\geq 0,..., y_r\geq 0;\ x_1\geq 0,..., x_{s-1} \geq 0;\
x_{11}\geq 0,...,x_{r(s-1)}\geq 0.
\end{cases}
\end{eqnarray}
For short, we set
$$\uda = (u_0,...,u_{rs+s-1}) = (z, y_1,...,y_r,
x_1,...,x_{s-1},x_{11},...,x_{r(s-1)}).$$
By Remark \ref{D2b}, a monomial
$\tda^{\bda} \in J \cap I^{n-1} $ if and only if the system (\ref{ED1}) has an integer
solution $\uda^*$ such that $u^*_0 = n,\ u^*_1= b_1,...,u^*_r = b_r$.

The corresponding system of homogeneous linear constraints is
\begin{eqnarray}\label{ED2}
\begin{cases}
y_j \geq a_{1j}x_1 + \cdots + a_{(s-1)j}x_{s-1} + a_{sj}(z-x_1-\cdots -
x_{s-1}), \ (j=1,...,r),\\
z \geq x_1+ \cdots + x_{s}, \\
y_j \geq a_{1j}x_{i1} + \cdots + a_{(s-1)j}x_{i(s-1)} + a_{sj}(z-x_{i1}-\cdots
- x_{i(s-1)}), \\ (i,j=1,...,r;\ j\neq i),\\
z \geq x_{i1}+ \cdots + x_{i(s-1)}, \ (i=1,...,r), \\
z \geq 0;\ y_1\geq 0,..., y_r\geq 0;\ x_1\geq 0,..., x_{s-1} \geq 0;\
x_{11}\geq 0,...,x_{r(s-1)}\geq 0.
\end{cases}
\end{eqnarray}
An integer  solution $(n, \bda, \xda)$ of this system gives a monomial $\tda^{\bda} \in
J \cap I^n = I^n$. Denote the sets of all integer solutions of (\ref{ED1}) and
(\ref{ED2}) by $E$ and $S$, respectively. Then $K[S] , K[E] \subseteq K[\uda]$
and $K[E]$ is a $K[S]$-module. Equip $K[S]$ and $K[E]$ with an $\Nset$-grading
by setting
$$\deg(\uda^{\cda}) = c_0.$$
Let $\Ical$ be the ideal of $K[S]$ generated by all binomials $\uda^{\alphda} -
\uda^{\betda}$,  such that $\alpha_0= \beta_0,...,\alpha_r = \beta_r$.

\begin{Lemma}\label{D2} There is an isomorphism of $\Nset$-graded rings
$$K[S]/\Ical \cong \Rcal := \oplus_{n\geq 0} I^nt^n.$$
\end{Lemma}

\begin{pf} The above discussion shows that there is an epimorphism of $\Nset$-graded
rings:
$$\begin{array}{rl}
K[S] & \rightarrow \Rcal,\\
\uda^{\cda} & \mapsto t_1^{c_1}\cdots t_r^{c_r}t^{c_0}.
\end{array}$$
The kernel of this map is exactly $\Ical$. The proof is similar to that of
Lemma 4.1 in \cite{St}, or we can argue directly as follows. By Lemma
\ref{A1}, $K[S]$ is generated by a finite number of monomials, say
$\uda^{\cda_1},..., \uda^{\cda_p}$. Consider the polynomial ring $K[\vda]$ of
$p$ new variables $\vda = (v_1,...,v_p)$. By \cite{St}, Lemma 4.1, the kernel
of the epimorphism
$$\psi:\ K[\vda]\rightarrow K[S],\ \ \psi(v_i) =
\uda^{\cda_i},$$
 is the ideal $I_{\Acal}$ generated by binomials
$\vda^{\alphda} - \vda^{\betda}$ such that $\sum_{i=1}^p \alpha_i \cda_i =
\sum_{i=1}^p \beta_i \cda_i$. Such an ideal is called {\it toric ideal}
associated to the matrix $\Acal : = \{ \cda_1,...,\cda_p\}$. Let $\cda'_i =
(c_{i0},...,c_{ir})$ and $\Acal' : = \{ \cda'_1,...,\cda'_p\}$. Again by
\cite{St}, Lemma 4.1, the kernel of the epimorphism
$$\chi: \ K[\vda]  \rightarrow \Rcal,\ \
\vda_i  \mapsto t_1^{c_{i1}}\cdots t_r^{c_{ir}}t^{c_{i0}},$$
is $I_{\Acal'}$.
Clearly $I_{\Acal} \subseteq I_{\Acal'}$, $\psi(I_{\Acal'}) = \Ical$. Hence
$\chi$ induces an isomorphism
$$\varphi:\ K[S]  \rightarrow \Rcal ,$$
such that $\text{Ker}\, \varphi = \Ical$ and $\varphi (\uda^{\cda_i}) = \chi
(v_i)$. This implies
$$ \varphi (\uda^{\cda}) = t_1^{c_1}\cdots t_r^{c_r}t^{c_0}$$
for all $\cda \in S$.
\end{pf}

By this isomorphism, we can consider the quotient module $K[E]/ \Ical K[E]$ as
a module over $\Rcal$. Of course, $H^0_{\mm G}(G) $ can be considered  as a module over
$\Rcal$, too.

\begin{Lemma}\label{D3} Let $r\geq 2$. Then there is an epimorphism of $\Nset$-graded modules
over
$\Rcal$
$$K[E]/\Ical K[E]\rightarrow  \oplus_{n\geq 1} \frac{J \cap I^{n-1}}{I^n} t^n = H^0_{\mm G}(G)
.$$
\end{Lemma}

\begin{pf} The set $M= \oplus_{n\geq 1}(J \cap I^{n-1}) t^n$ is a module over
$\Rcal$ and contains the ideal $I\Rcal$. The isomorphism $\varphi$ in the
proof of Lemma \ref{D2} induces a homomorphism
$$\begin{array}{rl}
K[E]/\Ical K[E] & \rightarrow M ,\\
\uda^{\cda} & \mapsto t_1^{c_1}\cdots t_r^{c_r}t^{c_0},
\end{array}$$
which is clearly surjective. Since $H^0_{\mm G}(G)\cong M/I\Rcal$, it is an
image of $K[E]/\Ical K[E]$.
\end{pf}

\begin{Proposition}\label{D4} Let $r\geq 2$ and $d$ be the maximal degree of the generators of
$I$, i.e. $d = \max_i (a_{i1} + \cdots + a_{ir})$. Then the $\Rcal$-module $H^0_{\mm G}(G)$ is
generated by homogeneous elements of degrees less than
$$B_1 := d(rs+s+d)(\sqrt{r})^{r+1}(\sqrt{2} d)^{(r+1)(s-1)}.$$
\end{Proposition}

\begin{pf} By Lemma \ref{D3}, it suffices to show that $K[E]$ is generated over $K[S]$ by
monomials of degrees less than $B_1$. The system (\ref{ED1}) has $rs+s$
variables. Denote by $\delta (x)$ the vector obtained from the coefficient
vector of a variable $x$ by deleting already known zero entries. For
simplicity we write it in the row form. Then
$$\delta(x_{ik}) = (a_{k1}- a_{s1},...,a_{k(i-1)}- a_{s(i-1)}, a_{k(i+1)}-
a_{s(i+1)},...,a_{kr}- a_{sr}, -1).$$ We have
$$\begin{array}{ll}
\|\delta(x_{ik})\|^2 & \leq 1+ (a_{k1}- a_{s1})^2 +\cdots + (a_{kr}- a_{sr})^2
\\ &\leq 1 + (a_{k1}^2 +\cdots + a_{kr}^2) + (a_{s1}^2 +\cdots + a_{sr}^2) \\
& \leq 1+ (a_{k1} +\cdots + a_{kr})^2 + (a_{s1} +\cdots + a_{sr})^2 - 2
\sum_{i<j} a_{si}a_{sj}\\
& < 2d^2 \ \ (\text{by\ the \ condition}\ (\sharp)).
\end{array}$$
Similarly,
$$\| \delta(x_{i})\|^2  < 2d^2.$$
For all $j=1,..,r$, $\delta (y_j) = (1,...,1)$ ($r$ entries $1$). Hence
$$\| \delta(y_j) \|^2 = r.$$
Further,
$$\delta(z) = (\ada_s,\ada_s[1],...,\ada_s[r], 1,...,1) \ \ (r+1 \
\text{entries} \ 1).$$
 This yields
$$\|\delta(z)\|^2 = r(a_{s1}^2 +\cdots + a_{sr}^2) + r+1 < rd^2.$$
For the free coefficients of (\ref{ED1}) we have $\delta =
(a_{s1},...,a_{sr},1)$. So
$$ \|\delta \|^2 < d^2.$$
Applying Lemma \ref{A2} we get that $K[E]$ is generated over $K[S]$ by
monomials $\uda^{\cda}$ with
$$\|\cda\|_*  < (rs+s+d) \sqrt{r}d (\sqrt{2}d)^{(r+1)(s-1)} \sqrt{r}^r = B_1.
$$
Since $\deg(\uda^{\cda}) = c_0 \leq \|\cda\|_* < B_1$, the proof of the
proposition is completed.
\end{pf}

\begin{Proposition}\label{D5} Keep the notation of Proposition \ref{D4}. Let
$n\geq B_1$ be an integer. Then
$$\Ass(I^n /I^{n+1}) \supseteq \Ass(I^{n+1}/I^{n+2}).$$
\end{Proposition}

\begin{pf} Induction on the number of variables $r$. The case $r=1$ is
trivial. Let $r\geq 2$.
  By the induction hypothesis we have
$$\cup_{i=1}^r \Ass(I[i]^n/I[i]^{n+1}) \supseteq \cup_{i=1}^r \Ass(I[i]^{n+1}/I[i]^{n+2}).$$
If $\mm \in \Ass(I^n /I^{n+1})$,  then the above inclusion together with Lemma
\ref{B1} obviously give the assertion. Let $\mm \not\in \Ass(I^n /I^{n+1})$.
Then $H^0_{\mm G}(G)_n= H^0_{\mm}(I^n /I^{n+1}) = 0$. Since the module
$H^0_{\mm G}(G)$ is generated by elements of degrees less than $B_1$ over the
standard graded ring $\Rcal$ and $n\geq B_1$,  we must have $H^0_{\mm
G}(G)_{n+1} = 0$. This implies $\mm \not\in \Ass(I^{n+1}/I^{n+2})$. Hence, we
have by Lemma \ref{B1}
$$\Ass(\frac{I^{n+1}}{I^{n+2}}) = \Ass(\frac{I^{n+1}}{I^{n+2}})\setminus \{\mm\} = 
\cup_{i=1}^r
\Ass(\frac{I[i]^{n+1}}{I[i]^{n+2}}) \subseteq \cup_{i=1}^r
\Ass(\frac{I[i]^n}{I[i]^{n+1}})\subseteq \Ass(\frac{I^n }{I^{n+1}}).$$
\end{pf}

In order to get the reverse inclusion we use a result of S. McAdam and P.
Eakin (see  \cite{ME}, pp. 71, 72 and also \cite{Sh}, Proposition 2.4). Let
$$\Rcal_+ = \oplus_{n>0}I^nt^n.$$
The local cohomology module $H^0_{\Rcal_+}(G)$ is also a $\Zset$-graded
$\Rcal$-module. Let
$$a_0(G) = \sup\{n | \ H^0_{\Rcal_+}(G)_n \neq 0 \}.$$
(This number is to be taken as $-\infty$ if $H^0_{\Rcal_+}(G)=0$.) It is
related to an important invariant called the Castelnuovo-Mumford regularity of $G$
(see, e.g., \cite{Sh}). We have

\begin{Lemma}{\rm (\cite{ME}, Proposition 2.4)}. \label{D6} $ \Ass(I^n /I^{n+1}) \subseteq
\Ass(I^{n+1}/I^{n+2})$ for all $n> a_0(G)$.
\end{Lemma}

To define $H^0_{\Rcal_+}(G)$, let us recall the Ratliff-Rush closure of an
ideal:
$$\widetilde{I^n} = \cup_{m\geq 1} I^{n+m}:I^m.$$
This immediately gives

\begin{Lemma} \label{D7} For all $n> 0$ we have $H^0_{\Rcal_+}(G)_{n-1}\cong
(\widetilde{I^n} \cap I^{n-1})/ I^n$.
\end{Lemma}

Recall that $I = (\tda^{\ada_1},...,\tda^{\ada_s})$.

\begin{Lemma} \label{D8} For all $n> 0$ we have
$$\widetilde{I^n} = \cup_{m\geq 0} I^{n+m}:
(\tda^{m\ada_1},...,\tda^{m\ada_s}).$$
\end{Lemma}

\begin{pf} Since $\tda^{m\ada_i} \in I^m$, the inclusion $\subseteq$ is
obvious. To show the inclusion $\supseteq$, let
$$x\in I^{n+m}:
(\tda^{m\ada_1},...,\tda^{m\ada_s}).$$ Put $m' = sm$ and let $y$  be an
arbitrary element in $I^{m'}$. Then $y= (\tda^{m\ada_i})y'$ for some $i$ and
$y' \in I^{m'-m}$. We have
$$xy = y'(x \tda^{m\ada_i}) \in y'I^{n+m} \subseteq I^{n+ m'}.$$
This implies
$$x \in I^{n+ m'}: I^{ m'} \subseteq \widetilde{I^n}.$$
\end{pf}

\begin{Proposition}\label{D9} We have
$$a_0(G) < B_2 := s(s+r)^4 s^{r+2} d^2(2d^2)^{s^2-s+1}.$$
\end{Proposition}

\begin{pf} Consider the following system of linear constraints
\begin{eqnarray}\label{ED3}
\begin{cases}
y_j + a_{ij} x \geq a_{1j}x_{i1} + \cdots + a_{(s-1)j}x_{i(s-1)} + a_{sj}(z+x
- x_{i1} -\cdots - x_{i(s-1)}), \\
z +x \geq x_{i1}+ \cdots + x_{i(s-1)}, \\
 (i=1,...,s;\ j=1,...,r),\\
z\geq 0,\ x \geq 0;\ y_1\geq 0,..., y_r\geq 0;\ x_{11}\geq
0,...,x_{s(s-1)}\geq 0.
\end{cases}
\end{eqnarray}
By Lemma \ref{D8} and Remark \ref{D2b} we have $\tda^{\bda} \in
\widetilde{I^n}$ if and only if there is an integer  solution
$$\uda := (u_0,...,u_{s(s-1)+r+1}) := (z,x,y_1,...,y_r,x_{11},...,x_{s(s-1)})$$
of (\ref{ED3}) such that $z=n$ and $\bda = (u_2,...,u_{r+1})$. This system has
$s(s-1)+r+2$ variables. Using the notation in the proof of Proposition
\ref{D4}, a straightforward calculation gives
$$ \|\delta(x_{ij}) \|^2 < 2d^2,\ \|\delta(y_k) \|^2 = s,\ \|\delta(z) \|^2 <
sd^2,\ \|\delta(x) \|^2 < 2sd^2.$$ Let $S$ be the set of all integer solutions of
(\ref{ED3}). By Lemma \ref{A1}, the ring $K[S]$ is generated by monomials, say
$\uda^{\cda_1},..., \uda^{\cda_p}$, with
$$\begin{array}{ll}
 \|\cda_j\|_* & < (s(s-1)+r+2)(\sqrt{2}d)^{s(s-1)}\sqrt{s}d \sqrt{2s}d
 \sqrt{s}^r\\
 & < [(s+r)^2 d(\sqrt{2}d)^{s^2-s+1}\sqrt{s}^{r+2}] - 1 =: B_3.
 \end{array}$$
Fix such a generator $\uda^{\cda_j}$ of $K[S]$. Let
$$\cda'_j = (c_{j2},...,c_{j(r+1)}).$$
 Then from (\ref{ED3}) we have
$$(\tda^{\ada_i})^{c_{j1}} \tda^{\cda'_j} \in I^{c_{j0} + c_{j1}}.$$
Since $c_{j1}< B_3$, this implies that
$$(\tda^{\ada_i})^{B_3} \tda^{\cda'_j} \in I^{c_{j0} + B_3},$$
for all $i\leq s$. Let
$$B_4 = sB_3.$$
Since
$$I^{sB_3} = \sum_{i=1}^s (\tda^{\ada_i})^{B_3}I^{(s-1)B_3},$$ from the
above relationship we get
\begin{equation}\label{ED4}
I^{B_4} \tda^{\cda'_j} \subseteq  I^{c_{j0} + B_4}.
\end{equation}
We will show that $\widetilde{I^n} = I^n$ for all $ n\geq B_4(B_3+1)$. For
this aim, let $\tda^{\bda} \in \widetilde{I^n}$ with
$$n\geq B_4(B_3+1).$$
 Then
there are  $\alpha \in \Nset$ and $\alphda = (\alpha_{11},...,
\alpha_{s(s-1)})$,  such that $(n,\alpha,\bda, \alphda) \in S$. Since
$\cda_1,...,\cda_p$ generate $S$, there are nonnegative integers $m_1,...,m_p$
such that
$$(n,\alpha,\bda,\alphda) = m_1\cda_1+ \cdots + m_p\cda_p.$$
This implies
$$m_1 c_{10}+ \cdots + m_p c_{p0} = n,$$
and
$$\tda^{\bda} = (\tda^{\cda'_1})^{m_1} \cdots (\tda^{\cda'_p})^{m_p}.$$
Repeated application of (\ref{ED4}) gives
$$(\tda^{\ada_i})^{B_4}\tda^{\bda} \in I^{n + B_4},$$
for all $i \leq s$. In other words,
$$\tda^{\bda} \in I^{n + B_4} : (\tda^{B_4\ada_1},...,\tda^{B_4\ada_s}).$$
Hence there is $\betda = (\beta_{ij})\in \Nset^{s(s-1)}$ such that $(n, B_4,
\bda, \betda) \in S$. Then one can write
\begin{equation}\label{ED5}
(n, B_4,\bda, \betda) = \sum_{c_{i1}=0} m'_i\cda_i + \sum_{c_{i1}> 0}
m''_i\cda_i,
\end{equation}
for some $m'_i, m''_i \in \Nset$. We have
$$n = \sum_{c_{i1}=0} m'_ic_{i0} + \sum_{c_{i1}> 0}
m''_ic_{i0}.$$
Comparing the second components in (\ref{ED5}) gives
$$ \sum m''_i \leq \sum m''_ic_{i1} = B_4.$$
Since $c_{i0} \leq \| \cda_i\|_* \leq B_3$, we must have
$$ \sum m'_ic_{i0}  = n -\sum_{c_{i1}> 0}m''_ic_{i0} \geq n- B_3 \sum
m''_i \geq B_4(B_3+1) - B_3B_4 = B_4.$$
 In particular,  the set of $\cda_i$ with
$c_{i1}=0$ in (\ref{ED5}) is not empty. From (\ref{ED3}) one can see that for
such an index $i$ we have $\tda^{\cda'_i} \in I^{c_{i0}}$. Therefore, by
(\ref{ED5}) one obtains
$$\tda^{\bda} = \tda^{\bda'}\tda^{\bda''},$$
where
$$\tda^{\bda'} \in I^{\sum m'_ic_{i0}} \subseteq I^{B_4},$$
and
$$\tda^{\bda''} = \prod_{c_{i1}>0}(\tda^{\cda'_i})^{m''_i}.$$
Repeated application of (\ref{ED4}) once more gives us $\tda^{\bda} \in I^n$.
Thus we have shown
$$\widetilde{I^n} = I^n.$$
 Since
$$B_4(B_3+1) = sB_3(B_3+1) < s [(s+r)^2
d(\sqrt{2}d)^{s^2-s+1}\sqrt{s}^{r+2}]^2 \leq B_2,$$ from Lemma \ref{D7} we get
$a_0(G) < B_2$.
\end{pf}

Finally we can prove

\begin{Theorem} \label{D10} Let
$$B = \max\{d(rs+s+d)(\sqrt{r})^{r+1}(\sqrt{2} d)^{(r+1)(s-1)}, \
s(s+r)^4 s^{r+2} d^2(2d^2)^{s^2-s+1} \}.$$
Then we have
$$\Ass(I^n/I^{n+1}) = \Ass(I^B/I^{B+1})$$
for all $n\geq B$.
\end{Theorem}

\begin{pf} Note that $B= \max\{B_1, B_2\}$.  By Proposition \ref{D5}, $\Ass(I^n/I^{n+1})
\subseteq
\Ass(I^B/I^{B+1})$. By Lemma \ref{D6} and Proposition \ref{D9},
$\Ass(I^n/I^{n+1}) \supseteq \Ass(I^B/I^{B+1})$.
\end{pf}

The number $B$  in the above theorem  is very big. However
the following examples  show that such a number $B$ should
depend
 on $d$ and $r$.

\begin{Example} \label{E11} {\rm Let $d\geq 4$ and
$$I = (x^d,x^{d-1}y, xy^{d-1}, y^d,x^2y^{d-2}z) \subset K[x,y,z].$$
A monomial of this type was used in the proof of Theorem 4.1 in \cite{HH}. We
have
$$I^n: (x,y,z)^\infty = (x^d,x^{d-1}y, xy^{d-1}, y^d,x^2y^{d-2})^n = I^n $$
if and only if $n\geq k-2$. Hence
$$\Ass(I^{n-1}/I^n) = \begin{cases} \{ (x,y,z),
(x,y)\} & \ \text{if} \ n < d-2, \\
 \{  (x,y)\} & \ \text{if} \ n \geq  d-2.
 \end{cases}$$ }
\end{Example}

\begin{Example} \label{E12} {\rm Let $r\geq 4$ and $d>r-3$.
We put
$$u = t_1^{r-3 \choose 0} t_2^{r-3 \choose 1}\cdots  t_{r-3}^{r-4 \choose 0} \ \text{and} \
v= t_1^{\beta_1}\cdots  t_{r-3}^{\beta_{r-4}}t_{r-2}^{d-r+2},$$
where
$$\beta_i =   \begin{cases}  0 & \text{if}\ r-3-i\  \text{is\ even},\\
2{r-3\choose i} & \text{if}\ r-3-i\  \text{is\ odd}.
 \end{cases}$$
Let
$$I= (ut_1^d,ut_2^{d-1}t_r,...,ut_{r-2}^{d-r+3}t_r^{r-3}, ut_{r-1}t_r^{d-1}, vt_r^{r-3}),$$
and  $J$ be the integral closure $\overline{I^r}$  of $I^r$.  Assume that
$$\Ass(J^{n-1}/J^n) = \Ass(J^{B-1}/I^B)$$
for all $n\geq B$. Then
$$B \geq  \frac{d(d-1)\cdots (d-r+3)}{r(r-3)}.$$
Note that in this example $J$ is generated by monomials of degree $r(d+2^{r-3} -1)$.  Thus, if $r$
is fixed, then $B$ is at least $O(d(J)^{r-2})$, where $d(J)$ is the maximal degree of the generators of
$J$.
}\end{Example}

\begin{pf}
By \cite{V}, Corollary 7.60, $J$ is a normal ideal. Using the filtration
$$ J^n =  \overline{I^{rn}} \subset  \overline{I^{rn -1 }} \subset \cdots \subset  \overline{I^{r(n-1)}} 
= J^{n-1} ,$$
we get
$$\Ass( \overline{I^{rn-1}}/ \overline{I^{rn}} ) \subseteq \Ass (J^{n-1}/J^n) \subseteq 
\cup_{i=1}^r \Ass( \overline{I^{r(n-1)+i-1}} / \overline{I^{r(n-1)+i}}) .$$
By virtue of \cite{T}, Proposition 4, it is shown in the proof of \cite{T}, Proposition 16, that
$$\mm \in \Ass( \overline{I^{k-1}}/ \overline{I^{k}}) \ \ \text{ for\  all}\ \ k\gg 0,$$
and
$$\mm \not\in \Ass( \overline{I^{k-1}}/ \overline{I^{k}}) \ \ \text{if}\ \ k< \delta := 
\frac{d(d-1)\cdots (d-r+3)}{(r-3)}.$$
Hence $\mm \in \Ass (J^{n-1}/J^n) $ for all $n\gg 0$.  Assume that
$$B < \frac{d(d-1)\cdots (d-r+3)}{r(r-3)} = \delta /r.$$
Then $Br < \delta $, and hence $\mm \not\in  \Ass( \overline{I^{r(B-1)+i-1}} / 
\overline{I^{r(B-1)+i}})$ for all $i\leq r$.  This implies $\mm\not\in  \Ass (J^{B-1}/J^B)$, a 
contradiction.  \end{pf}

\end{document}